\documentclass[12pt,twoside]{amsart}
\usepackage{amssymb}
\usepackage{amscd}

\title[$4$-fold semi-stable log flips]
{On Termination of 4-fold semi-stable log flips}  
\author{Osamu Fujino} 
\address{Graduate School of Mathematics\\ 
 Nagoya University, Chikusa-ku Nagoya 464-8602 Japan}
\email{fujino@math.nagoya-u.ac.jp}
\newcommand{\xdiscrep}[0]{{\operatorname{discrep}}}
\newcommand{\xSupp}[0]{{\operatorname{Supp}}}
\newcommand{\Exc}[0]{{\operatorname{Exc}}}
\newtheorem{thm}{Theorem}[section]
\newtheorem{lem}[thm]{Lemma}

\newtheorem{prop}[thm]{Proposition}
\newtheorem{lem-def}[thm]{Definition-Lemma}
\newtheorem{assu}[thm]{Assumption}

\theoremstyle{definition}
\newtheorem{ex}[thm]{Example}
\newtheorem{defn}[thm]{Definition}
\newtheorem{rem}[thm]{Remark}
\newtheorem*{ack}{Acknowledgments}     
\newtheorem*{notation}{Notation}         
       
\newtheorem*{case}{Case}         
       
\newtheorem{step}{Step}
\newtheorem{say}[thm]{}
\begin{document}
\bibliographystyle{amsalpha+}

\begin{abstract}
In this paper, we prove the termination of 
$4$-fold semi-stable log flips under the assumption 
that there always exist $4$-fold (semi-stable) log 
flips. 
\end{abstract}

\maketitle

\section{Introduction}\label{sec1}

One of the most important conjectures in the (log) minimal 
model program (MMP, for short) is {\em{$($log$)$ Flip 
Conjecture}} II. It claims that 
any sequence of (log) flips: 
$$
\begin{matrix} 
(X_0,B_0) & \dashrightarrow & {(X_1,B_1)} & 
\dashrightarrow & {(X_2,B_2)} 
&\dashrightarrow\cdots\\
{\ \ \ \ \searrow} & \ &  {\swarrow}\ \  {\searrow} & \ &  
{\swarrow\ \ } &\\
 \ & Z_0 & \  \ & Z_1 & \ & &, 
\end{matrix}
$$
has to terminate after 
finitely many steps. For the non-log case, the conjecture in 
dimension $4$ was proved for the terminal flips by Kawamata 
in \cite{kmm},  
and for the flops by Matsuki in \cite{matsuki}. 
For the log case, we proved it for 
$4$-fold canonical flips in \cite{fujino}, 
which is a first step 
to prove the log Flip Conjecture II in dimension $4$. 
We note that the main theorem of \cite{fujino} 
contains the above mentioned 
results of Kawamata and Matsuki. 

Recently, Shokurov treats the log Flip Conjecture II 
in a much more general setting. For the details, 
see \cite{S2} and \cite{S3}. 
 
The main purpose of this paper is to prove the 
following theorem, which is a $4$-dimensional analogue of 
\cite[Theorem 7.7]{km}, under the assumption 
that there always exist $4$-fold (semi-stable) log 
flips (see Assumption \ref{sho} below). 
We will prove it 
by the {\em{crepant descent}} technique by Kawamata and 
Koll\'ar. 
For the details of the (log) semi-stable MMP, 
see \cite[\S 7.1]{km}. 

\vspace{3mm}

We will work over $\mathbb C$, the complex number field, 
throughout this paper. 

\begin{thm}[Termination of $4$-fold semi-stable log 
flips]\label{main}
Let $(X,B)$ be a $\mathbb Q$-factorial 
projective $4$-dimensional dlt pair, 
$f:X \longrightarrow Y$ a projective surjective 
morphism and 
$g:Y\longrightarrow C$ a flat morphism 
to a non-singular curve $C$ such that 
$h:=g\circ f :(X,B)\longrightarrow C$ is a dlt morphism 
$($for the definition of 
dlt morphisms, see {\em{Definition \ref{plt-mo}}} 
below$)$. 
Then an arbitrary sequence of extremal 
$(K_X+B)$-flips over $Y$ is finite. 
\end{thm}

In the proof of Theorem \ref{main}, we need the 
following assumption: Assumption \ref{sho}. We note that 
all the flips we need here is $4$-fold semi-stable (log) flips, 
which is a special case of Assumption \ref{sho} (see Definition 
\ref{semi} and Remark \ref{exist} below).  
Recently, Shokurov announced a proof of 
the existence of $4$-fold log flips in \cite{S}. 
So, this assumption seems to be reasonable. 
We recommend the readers to see \cite{S}.  

\begin{assu}\label{sho}
Let $(X,B)$ be a $4$-dimensional klt pair and 
$f:X\longrightarrow Z$ a flipping contraction with 
respect to $K_X+B$. Then 
$f$ has a flip. 
\end{assu}

For the proof of Theorem \ref{main}, we need the 
following two theorems. First, 
we recall the special termination theorem. 
For the details, see \cite[Section 2]{S} and \cite{fujino1}. 

\begin{thm}[$4$-dimensional special termination]\label{st}
Let $(X,B)$ be a $\mathbb Q$-factorial dlt $4$-fold. 
Consider a sequence of extremal $(K_{X_i}+B_i)$-flips 
starting from $(X,B)=(X_0,B_0)${\em{:}}
$$
\begin{matrix} 
(X_0,B_0) & \dashrightarrow & {(X_1,B_1)} & 
\dashrightarrow & {(X_2,B_2)} 
&\dashrightarrow\cdots\\
{\ \ \ \ \searrow} & \ &  {\swarrow}\ \  {\searrow} & \ &  
{\swarrow\ \ } &\\
 \ & Z_0 & \  \ & Z_1 & \ & &, 
\end{matrix}
$$ 
Then after finitely many flips, flipping locus {\em{(}}and 
thus the flipped locus{\em{)}} is disjoint from $\llcorner B_i\lrcorner$. 
\end{thm}

Next, the following theorem is the main theorem of \cite{fujino}. 

\begin{thm}[Termination of $4$-fold canonical flips]\label{can}
Let $X$ be a normal projective $4$-fold and 
$B$ an effective $\mathbb Q$-divisor such 
that $(X,B)$ is canonical, that is, $\xdiscrep (X,B)\geq 0$. 
Consider a sequence of $(K_{X_i}+B_i)$-flips starting from 
$(X,B)=(X_0,B_0)${\em{:}}
$$
\begin{matrix} 
(X_0,B_0) & \dashrightarrow & {(X_1,B_1)} & 
\dashrightarrow & {(X_2,B_2)} 
&\dashrightarrow\cdots\\
{\ \ \ \ \searrow} & \ &  {\swarrow}\ \  {\searrow} & \ &  
{\swarrow\ \ } &\\
 \ & Z_0 & \  \ & Z_1 & \ & &, 
\end{matrix}
$$
Then this sequence terminates after finitely many steps. 
\end{thm} 

We note that we don't need Assumption 
\ref{sho} in the proofs of 
Theorems \ref{st} and \ref{can}. 

\begin{ack} 
I would like to thank Professor V.~V.~Shokurov for 
comments. 
\end{ack}

\begin{notation}
Let $\mathbb Z_{>0}$ (resp.~$\mathbb Z_{\geq 0}$) be 
a set of positive (resp.~non-negative) integers. 
For $d\in \mathbb Q$, let 
$\llcorner d \lrcorner =\max \{t \in \mathbb Z 
\ | \ t\leq d\}$ and $\{ d \}=d-\llcorner d\lrcorner$. 
Let $D=\sum d_iD_i$ be a $\mathbb Q$-divisor such 
that all the $D_i$'s are distinct. 
We put $\llcorner D\lrcorner =\llcorner d_i\lrcorner D_i$ 
(the {\em{round down} of $D$}) 
and $\{D\}=\sum \{d_i\}D_i$ (the {\em{fractional part of}} 
$D$). 
\end{notation}

\section{Preliminaries}\label{sec2}

In this section, we collect basic properties and 
definitions. 

\begin{say} First, let us recall the definitions of 
discrepancies and singularities of pairs. 

\begin{defn}[Discrepancies and singularities for pairs]\label{def1}
Let $X$ be a normal variety and $D=\sum d_i D_i$ a 
$\mathbb Q$-divisor on $X$, 
where $D_i$ is irreducible for 
every $i$ and $D_i\ne D_j$ for 
$i\ne j$, such that $K_X+D$ is $\mathbb 
Q$-Cartier. 
Let $f:Y\longrightarrow X$ be a proper birational morphism 
from a normal variety $Y$. 
Then we can write 
$$
K_Y=f^{*}(K_X+D)+\sum a(E,X,D)E, 
$$ 
where the sum runs over all the distinct prime divisors $E\subset Y$, 
and $a(E,X,D)\in \mathbb Q$. This $a(E,X,D)$ is called the 
{\em discrepancy} of $E$ with respect to $(X,D)$. 
We define 
$$
\xdiscrep (X,D):=\inf _{E}\{a(E,X,D)\ 
|\  E \text{ is exceptional over}\  X  \}.
$$ 
From now on, we assume 
that $0\leq d_i \leq 1$ 
for every $i$. 
We say that $(X,D)$ is 
$$
\begin{cases}
\text{terminal}\\
\text{canonical}\\
\text{klt}\\
\text{plt}\\
\text{lc}\\
\end{cases}
\quad {\text{if}} \quad \xdiscrep (X,D) 
 \quad
\begin{cases}
>0,\\
\geq 0,\\
>-1\quad {\text {and \quad $\llcorner D\lrcorner =0$,}}\\
>-1,\\
\geq -1.\\
\end{cases}
$$ 
Here klt is short for {\em {Kawamata log terminal}}, plt 
for {\em{purely log terminal}}, and lc for 
{\em{log canonical}}. 

If there exists a log resolution $f:Y\longrightarrow X$ 
of $(X,D)$, that is, $Y$ is non-singular, the exceptional 
locus $\Exc(f)$ is a divisor, and $\Exc(f)\cup 
f^{-1}(\xSupp D)$ is a simple normal crossing 
divisor, such 
that $a(E_i,X,D)>-1$ for every exceptional divisor $E_i$ on 
$Y$, then the pair $(X,D)$ is called {\em{dlt}}. 
Here, dlt is short for {\em{divisorial log terminal}}. 
\end{defn}
\end{say}

\begin{say} 
Next, let us recall the definition of dlt morphisms and 
define plt morphisms. 

\begin{defn}[{\cite[Definition 7.1]{km}}]\label{plt-mo}
Let $X$ be a normal variety, $B$ an effective $\mathbb Q$-divisor 
on $X$ and $f:X \longrightarrow C$ a non-constant morphism to 
a non-singular curve $C$. 
We say that $f:(X,B)\longrightarrow C$ is dlt (resp.~plt) 
if $(X,B+f^{*}P)$ is dlt (resp.~plt) for every closed point $P\in 
C$. 
We note that if $(X,B)\longrightarrow C$ is plt, then 
$(X,B)$ is klt. 
\end{defn}

The following lemma is a variant of adjunction and 
the inversion of 
adjunction. 
For the proof, see \cite[Theorem 5.50 (1), Proposition 5.51]{km}. 

\begin{lem}\label{pmor}
Let $(X,B)$ be a klt pair and $f:(X,B)\longrightarrow C$ 
a dlt morphism. 
Then the following four conditions are equivalent. 
\begin{itemize}
\item[(1)] $f:(X,B)\longrightarrow C$ is a plt morphism. 
\item[(2)] every connected component of any fiber is irreducible. 
\item[(3)] $(F, B|_F)$ is a klt pair for any fiber $F$. 
\item[(4)] all the fibers of $f$ are normal. 
\end{itemize}
\end{lem}

The next lemma is an analogue of \cite[Lemma 7.2 (4)]{km}. 
It easily follows from the definition of dlt pairs 
(see \cite[Definition 2.37]{km}). 
We leave the details to the readers. 

\begin{lem}\label{new}
Let $(X,B)$ be a klt pair and $f:(X,B)\longrightarrow C$ 
a dlt morphism. If $E$ is an exceptional divisor over $X$ 
such that the center of $E$ on $X$ is contained in 
a fiber, then the discrepancy $a(E,X,B)>0$. 
\end{lem}

We note the following properties, which is an easy consequence of 
the negativity lemma (cf.~\cite[Lemma 3.38]{km}). 

\begin{lem}[{cf.~\cite[Corollary 3.44]{km}}]\label{nega}
Let $\phi:(X,B)\dashrightarrow (X^+,B^+)$ be either 
a $(K_X+B)$-flip over $Y$ or 
a divisorial contraction of a $(K_X+B)$-negative 
extremal ray over $Y$, 
$f:Y\longrightarrow C$ a flat morphism onto 
a non-singular curve $C$, and 
$h:
=f\circ g:(X,B)\longrightarrow C$ is a dlt $($resp.~plt$)$ morphism. 
Then so is $h^+:(X^+,B^+)\longrightarrow C$. 
\end{lem}
\end{say}

\begin{say} 
Finally, we define {\em{semi-stable log 
flips}} (cf.~\cite[Theorem 7.8]{km}). 

\begin{defn}\label{semi}
Let $(X,B)$ be a $\mathbb Q$-factorial 
dlt pair and $f:X\longrightarrow W$ 
an {\em{extremal}} flipping contraction with 
respect to $K_X+B$. 
Here, ``extremal'' means that $X$ is $\mathbb Q$-factorial and 
the relative Picard number 
$\rho (X/W)=1$. 
Assume that there exists a flat morphism 
$g:W\longrightarrow C$ to a smooth curve such that 
$h:=g\circ f$ is dlt. 
Then the flip $f^+:X^+\longrightarrow W$ of $f$: 
$$
\begin{matrix}
{X\ \ \ \ \ } & \dashrightarrow & {\ \ \ \ \ \ X^+} \\
{\ \ \searrow} & \ &  {\swarrow\ \ } \\
 \ & W & \ & 
\end{matrix}
$$ 
that is, 
\begin{itemize}
\item[(i)] $f^+$ is small, 
\item[(ii)] $K_{X^+}+B^+$ is $f^+$-ample, where 
$B^+$ is the strict transform of $B$, 
\end{itemize}
is called a {\em{semi-stable $($log$)$ flip}} of $f$. 
Furthermore, if $(X,B)$ is terminal, 
that is, $\xdiscrep (X,B)>0$, 
then we call $f^+$ a {\em{semi-stable terminal flip}} of $f$. 
\end{defn}
\end{say}

We treat only one example here. 

\begin{ex}[$4$-fold semi-stable flip]\label{rei}
Let $V$ be a projective $3$-fold with $\mathbb Q$-factorial 
terminal singularities and 
$$
\begin{matrix}
{V\ \ \ \ \ } & \dashrightarrow & {\ \ \ \ \ \ V^+} \\
{\ \ \searrow} & \ &  {\swarrow\ \ } \\
 \ & Z & \ & 
\end{matrix}$$ 
an extremal $K_V$-flip. 
We define $X:=V\times \mathbb P^1$, $X^+:=V^+\times 
\mathbb P^1$, and $W:=Z\times \mathbb P^1$. 
We put $Y:=C:=\mathbb P^1$. 
Then 
$$
\begin{matrix}
{X\ \ \ \ \ } & \dashrightarrow & {\ \ \ \ \ \ X^+} \\
{\ \ \searrow} & \ &  {\swarrow\ \ } \\
 \ & W & \ & 
\end{matrix}$$ 
is an extremal $4$-fold semi-stable flip over $Y$. 
We note that the second projection $X\longrightarrow C$ 
is a plt morphism. 
It is not difficult to see that $\rho(X/W)=1$ and 
$X$ is $\mathbb Q$-factorial. 
In this case, the flipping and flipped loci are 
dominant onto $C$. 
\end{ex}

\section{Preparation}\label{sec-p}

This section is a preparation of the proof of the main theorem: 
Theorem \ref{main}. 

\begin{say} 
We write a sequence of $4$-fold semi-stable flips over $Y$ 
as follows: 
$$
\begin{matrix}
(X,B)&=:
(X_0,B_0) & \dashrightarrow & {(X_1,B_1)} & \dashrightarrow & {(X_2,B_2)} 
&\dashrightarrow\cdots\\
&{\ \ \ \ \searrow} & \ &  {\swarrow}\ \  {\searrow} & \ &  
{\swarrow\ \ } &\\
& \ & W_0 & \  \ & W_1 & \ & &, 
\end{matrix}
$$
where $\phi_i: 
X_i\longrightarrow W_i$ is an extremal flipping contraction 
with respect to $K_{X_i}+B_i$ over $Y$ and 
$\phi_i^{+}:X_{i+1}\longrightarrow W_i$ is the flip of $\phi_i$ 
for every $i$. 

By the special termination theorem:~Theorem \ref{st}, 
all the flipping and flipped loci are disjoint 
from $\llcorner B_i\lrcorner$ after finitely many flips. 
Therefore, we can assume that 
all the flipping and flipped loci are disjoint from 
$\llcorner B_i\lrcorner$ for every $i$ by shifting the index 
$i$. 
So, we can replace $B_i$ with its fractional part 
$\{B_i\}$ 
and assume that $(X_i, B_i)$ is klt. 
From now on, we assume that $(X_i, B_i)$ is klt for every $i$. 
\end{say}

Let us recall the following definition. 

\begin{defn}[{\cite[6.6 Definition]{FA}}] 
Let $(X,B)$ be a klt $n$-fold. 
By \cite[Proposition 2.3.6]{km}, 
there are only finitely many exceptional divisors with 
non-positive discrepancies. 
The number of these divisors is denoted by $e(X,B)$. 
Thus $(X,B)$ is terminal if and only if 
$e(X,B)=0$ by the definition of terminal pairs. 
\end{defn}

\begin{say}
We prove Theorem \ref{main} by induction on $e(X,B)$. 

If $e(X,B)=0$, then $(X,B)$ is terminal. 
Thus a sequence of flips always terminates by 
Theorem \ref{can}. Therefore, we assume that 
the theorem holds for $e(X,B)\leq e-1$, 
and prove it for $e(X,B)=e$ case. 
We note that $e(X_i, B_i)\geq e(X_{i+1}, B_{i+1})$ 
for all $i$ by the negativity lemma 
(cf.~\cite[Lemma 3.38]{km}). 
\end{say}

\begin{say}\label{2.8}
First, we add $f^*P$ to $B$, where $P$ is a closed 
point of $C$. 
By Theorem \ref{st}, 
we can assume that all the flipping and flipped loci 
are not dominant onto $C$ after finitely many flips. 
Thus, by shifting the index $i$ we can assume that 
all the flipping and flipped loci are contained in some 
fibers. 
\end{say}

So, we can assume that there are no semi-stable flips like 
Example \ref{rei}. 

\begin{say}\label{2.9}
Next, we add $\sum _Pf^*P$ to $B$, 
where $P$ runs through all the closed points of 
$C$ such that $f^*P$ is not normal. 
By Theorem \ref{st} again, 
we can assume that all the flipping and 
flipped loci are disjoint from non-normal fibers. 
We note that the normality of fibers are preserved by flips 
(see Lemmas \ref{pmor} and \ref{nega}). 
Therefore, we can assume that there exists a non-empty 
Zariski open set $U$ of $C$ such that 
all the flips occur over this open set $U$ and 
$(X_i,B_i)\longrightarrow C$ is a plt morphism 
over $U$ (see Definition \ref{plt-mo}). 
\end{say}

We recall the definition of $r(X,B)$. 

\begin{defn}[{\cite[6.9.8 Definition]{FA}}]\label{aaru}
Let $(X,B)$ be a klt $n$-fold. 
We put 
$$
s(X,B):=\min \{a(E,X,B)>0 \ |\ E {\text{\ is exceptional over\ } X}\}. 
$$
Then we define 
$$
r(X,B):=(4\ulcorner s(X,B)^{-1}\urcorner)!\in \mathbb Z_{>0}. 
$$ 
\end{defn}

We generalize the invariants $e(X,B)$, $r(X,B)$, and 
$\xdiscrep(X,B)$ for plt morphisms. 
By Lemma \ref{pmor} (3), a plt morphism is a family of klt 
pairs. So, the following definition is natural. 

\begin{lem-def}\label{relative}
Let $f:(X,B)\longrightarrow C$ be a plt morphism. 
Then 
$$
0\leq e(f;(X,B)):=\max_Fe(F,B|_F)<\infty, 
$$
$$
r(f;(X,B)):=\max_Fr(F,B|_F)\in \mathbb Z_{>0}, {\text{\ \ and}}
$$ 
$$
-1<\xdiscrep (f;(X,B)):=\min _F \xdiscrep (F, B|F)\leq 1, 
$$ 
where $F$ runs through all the fibers of $f$. 
We note that $K_F+B|_F:=(K_X+B+F)|_F$ is klt by adjunction 
$($see {\em{Lemma \ref{pmor}}}$)$. 
\end{lem-def}

\begin{proof} 
Take a log resolution $g:Z\longrightarrow X$ of the pair $(X,B)$ 
as in \cite[Proposition 2.36 (1)]{km}. 
We write 
$$
K_Z+D-E=g^*(K_X+B), 
$$ 
where $D=\sum a_iD_i$ and $E=\sum b_jE_j$ 
are both effective and have no common 
irreducible component. Let $G=\sum G_k$ be the 
$g$-exceptional divisor such that $a(G_k,X,B)=0$ for 
every $k$.  
We can assume that $\xSupp (D\cup G)$ is non-singular. 
There exists a non-empty Zariski open set $U$ such that 
$f\circ g$ is smooth and 
$\xSupp (D\cup E\cup G)$ is relatively normal crossing over $U$. 
We can assume that $g(D_i)\longrightarrow C$, 
$g(E_j)\longrightarrow C$, and 
$g(G_k)\longrightarrow C $ are flat 
over $U$ for every $i$, $j$, and $k$ after shrinking $U$. 
Over this open set $U$, $e(F,B|_F)$, 
$r(F,B|_F)$ (more precisely, $s(F,B|_F)$), and 
$\xdiscrep (F,B|_F)$ 
are 
constant. Therefore, 
$ e(f;(X,B))$, $r(f;(X,B))$, 
and $\xdiscrep (f; (X,B))$ are well-defined and 
have the required properties. 
\end{proof}

The next proposition will play crucial roles in 
the proof of the main theorem. 

\begin{prop}\label{taisetsu}
Let $f:(X,B)\longrightarrow C$ be a plt morphism 
and $D$ a $\mathbb Q$-Cartier Weil divisor on $X$. 
Then $mD$ is Cartier if and only if 
so is $mD|_F$ for every fiber $F$. 
In particular, if $K_X$ is $\mathbb Q$-Cartier, 
then $mK_X$ is Cartier if and only if 
so is $mK_F$ for every fiber $F$. 
\end{prop}

\begin{proof}
See, for example, \cite[Lemma 2.1.7]{huy}. 
We note that $(X,B+F)$ is plt and $F$ is Cartier. Thus, 
in a neighborhood of $F$, $X$ is smooth in codimension two. 
So, $\mathcal O_X(mD)|_{F}\simeq \mathcal O_F(mD|_F)$ and 
$\mathcal O_X(m(K_X+F))|_F\simeq 
\mathcal O_F(mK_F)$ for every $m\in \mathbb Z_{\geq 0}$ 
(cf.~\cite
[Proposition 5.26]{km}). 
\end{proof}

We recall the result in \cite[6.11 Theorem]{FA}. 
For the proof, see \cite[(6.11.5)]{FA}. 

\begin{thm}\label{key}
Let $(V,\Delta)$ be a klt $3$-fold and 
$E$ a $\mathbb Q$-Cartier Weil divisor on $V$. 
Then $mE$ is Cartier for some 
$$
1\leq m
\leq r(V,\Delta)^{2^{e(V,\Delta)}}
\genfrac{(}{)}{}{}{3}{1+\xdiscrep (V,\Delta)}^{2^{e(V,\Delta)}-1}. 
$$
\end{thm}

\begin{thm}\label{miso}
Let $(X,B)$ be a klt $4$-fold and $f:(X,B)\longrightarrow 
C$ a plt morphism. Let $E$ be a $\mathbb Q$-Cartier Weil 
divisor on $X$. Then $ME$ is Cartier for 
$$
M:=M(f;(X,B)):= (\ulcorner\varphi(f;(X,B)) \urcorner)!\in \mathbb Z_{>0} ,  
$$
where 
\begin{multline*}
\varphi(f;(X,B))\\:=r(f;(X,B))^{2^{e(f;(X,B))}}\genfrac{(}{)}
{}{}{3}{1+\xdiscrep (f;(X,B))}^{2^{e(f;(X,B))}-1}. 
\end{multline*}
Let $U$ be a non-empty Zariski open subset of $C$. 
Then it is obvious that 
the restriction 
$f|_{f^{-1}(U)}:(X,B)|_{f^{-1}(U)}\longrightarrow U$ is 
a plt morphism and $M(f|_{f^{-1}(U)}; 
(X,B)|_{f^{-1}(U)})$ divides $M(f;(X,B))$. 
\end{thm}
\begin{proof}
It is obvious by Theorem \ref{key}. 
We note that if $E$ is not dominant onto $C$, 
then it is obvious that $E$ is Cartier. 
The latter statement is obvious by the definition of $M$. 
\end{proof}

\section{Proof of the main theorem}\label{sec3}

We go back to the proof of the main theorem:~Theorem \ref{main}. 
Our proof is similar to the proof of \cite[6.11 Theorem]{FA}. 

\begin{proof}[Proof of {\em{Theorem \ref{main}}}]
We start the proof of the main theorem. 
\begin{step}
First, we take a log resolution of $(X,B)$. 
We write $p:Z\longrightarrow X$ and 
$$
K_Z+p^{-1}_*B=p^{*}(K_X+B)+E-F,  
$$ 
where $E$ and 
$$F:=\sum _{a_i\geq 0}a_iF_i
$$ 
are effective exceptional divisors and 
have no common irreducible components. 
If necessary, we further blow up $Z$. 
Then we can assume that 
$\sum_{a_i\geq 0} F_i$ contains 
all the exceptional divisors whose discrepancies are 
non-positive, $\xSupp 
(p^{-1}_*B\cup \sum F_i)$ is smooth and  
$\xSupp 
(p^{-1}_*B\cup \sum F_i\cup (f\circ g\circ p)^{*}P)$ 
is simple normal crossing 
for every $P\in C$. 
We note that $F_i$ is dominant onto $C$ for every $i$ by 
Lemma \ref{new}. 
We can assume that $\sum F_i\ne 0$, that is, $e=e(X,B)>0$. 
We consider 
$$
f\circ g\circ p: 
(Z,D^{\varepsilon}):=
(Z,p^{-1}_*B+F+\varepsilon \sum _{i\ne 0}F_i)\longrightarrow C. 
$$ 
It is easy to check that $(Z,D^{\varepsilon})$ 
is terminal and $f\circ g\circ p: 
(Z,D^{\varepsilon})\longrightarrow C
$ is a dlt morphism for $0<\varepsilon 
\ll 1$. 
Run the log MMP over $X$. 
Then we obtain a sequence of flips and divisorial 
contractions over $X$: 
$$
(Z,D^{\varepsilon}):=(Z_0,D^{\varepsilon}_0)\dashrightarrow 
(Z_1,D^{\varepsilon}_1)\dashrightarrow \cdots 
\dashrightarrow (Z_k,D^{\varepsilon}_k)\dashrightarrow\cdots.    
$$
By Theorem \ref{sho}, flips exist and 
any sequence of flips terminates since 
$e(Z_k,D^{\varepsilon}_k)<e=e(X,B)$ for every $k$. 
Then we obtain a relative log minimal 
model $q:(Z',B')\longrightarrow 
X$, which satisfies the following conditions:  
\begin{itemize}
\item[(1)] $f\circ g\circ q:(Z',B')\longrightarrow C$ is a dlt morphism. 
\item[(2)] $f\circ g\circ q:(Z',B')\longrightarrow C$ is a 
plt morphism over $U$ (see \ref{2.9}). 
\item[(3)] $e(Z',B')=e(X,B)-1$. 
\item[(4)] $(Z',B')$ is a $\mathbb Q$-factorial klt pair. 
\item[(5)] $K_{Z'}+B'=q^{*}(K_X+B)$, that is, $q$ is a log crepant 
morphism. 
\item[(6)] the relative Picard numbers $\rho (Z'/X)=1$ 
and $\rho (Z'/W_0)=2$. 
\end{itemize}
We note that $\alpha:Z\dashrightarrow Z'$ 
is an isomorphism at 
the generic point of $F_0$ and 
contracts $E 
+\sum_{i\ne 0} F_i$. 
\end{step}
\begin{step}
We put $p_0:(Z^0_0,B^0_0):=(Z',B')\longrightarrow 
X=:X_0$. 
We assume that we already have 
$p_i:(Z^0_i,B^0_i):=(Z',B')\longrightarrow 
X_i$. 
Run the log MMP to $(Z^0_i,B^0_i)$ over $W_i$. 
We obtain a sequence of flips and divisorial contractions 
over $W_i$: 
$$
Z_i^{0}\dashrightarrow Z_i^{1}\dashrightarrow \cdots 
\dashrightarrow Z_i^{k_i},   
$$ 
and a 
log minimal model $(Z_i^{k_i}, B_i^{k_i})$ over  $W_i$. 
This is a so-called {\em{2 ray games}}.  
Since $(X_{i+1}, B_{i+1})$ is the log canonical 
model of $(Z^0_i,B^0_i)$ over 
$W_i$, there exists a morphism $q_{i}:Z_i^{k_i}\longrightarrow X_{i+1}$. 

\begin{case}[A]
If all the steps in the above log MMP are flips, 
then we have $K_{Z_i^{k_i}}+B_i^{k_i}=q^*_{i}(K_{X_{i+1}}+B_{i+1})$. 
We define $p_{i+1}:(Z^0_{i+1}, B^{0}_{i+1})
:=(Z_{i}^{k_i}, B_{i}^{k_i})\longrightarrow X_{i+1}$. 
We put $c_{i+1}=0$ in this case.  
\end{case}

\begin{case}[B] 
If a divisorial contraction occurs in the above log MMP, then 
it is not difficult to 
see that the final step 
$\beta:Z_{i}^{k_i-1}\dashrightarrow Z_{i}^{k_i}$ is a divisorial 
contraction and 
$q_{i}:Z_{i}^{k_i}\longrightarrow X_{i+1}$ is an isomorphism 
(cf.~\cite[Lemma 6.39]{km} and \cite[6.5.5 Proposition]{FA}). 
We note that other steps in the above log MMP are all flips. 
We also note that 
$$
K_{Z_i^{k_i-1}}+B_i^{k_i-1}=(q_i\circ\beta)^*(K_{X_{i+1}}+B_{i+1})+
c_{i+1}F_0 
$$ 
for $c_{i+1}>0$, where $F_0$ is the proper transform of 
$F_0$ on $Z_{i}^{k_i-1}$. 
Then we put $p_{i+1}:(Z_{i+1}^{0}, B_{i+1}^{0}): 
=(Z_{i}^{k_i-1}, B_{i}^{k_i-1}-c_{i+1}F_0 ) 
\longrightarrow X_{i+1}$. 
\end{case}
\end{step}
\begin{step}
If Case (B) occurs only finitely many times, then we 
can assume that all the steps are Case (A). 
Then we obtain an infinite sequence of flips with respect to 
$K_{Z^j_i}+B^j_i$. 
Since $e(Z^j_i,B^j_i)<e(X,B)$, it is impossible. 
So, Case (B) occurs infinitely many times. 
The coefficient of $F_0$, where $F_0$ is the 
proper transform of $F_0$ on $Z^0_{i+1}$, 
in $B_{i+1}^0$ is 
$$
a_0-\sum _{0\leq j\leq i}c_{j+1},  
$$ 
where $a_0:=-a(F_0, X,B)\geq 0$, that is, $a(F_0, X,B)\leq 0$. 
Let $U_{i+1}$ be a non-empty Zariski open set of $U$ such 
that flips $(X_j,B_j)\dashrightarrow (X_{j+1}, B_{j+1})$ 
occur over $U\setminus U_{i+1}$ for $0\leq j\leq i$. 
We note that it is sufficient to consider the 
coefficient of $F_0$ over $U_{i+1}$ since $F_0$ is irreducible 
and dominant onto $C$. 
Let $N$ be a positive integer such that $NB_0$ is a Weil 
divisor. 
Then $NB_i$ is also a Weil divisor for every $i$. 
By Theorem \ref{miso}, 
$MN(K_{X_i}+B_i)$ is a Cartier divisor over $U_i$ for every $i$, 
where $M:=M(h|_{(h)^{-1}(U)}; 
(X,B)|_{(h)^{-1}(U)})$. 
We note that $M(h|_{(h)^{-1}(U_{i})}; 
(X,B)|_{(h)^{-1}(U_i)})$ divides $M$ by Lemma \ref{nega} 
and that $(X_{i}, B_i)$ is isomorphic to $(X,B)$ over $U_i$. 
Thus $MNB^0_i$ is a Weil divisor over $U_i$ for every $i$. 
So, we have that $MNc_j\in \mathbb Z_{\geq 0}$ for every $j$. 
Therefore, after finitely many steps, 
the coefficient of $F_0$ in $B_{i+1}^0$ is negative, 
that is, the discrepancy $a(F_0, X_{i+1}, B_{i+1})>0$. 
Thus, $e(X_{i+1}, B_{i+1})<e=e(X,B)$. 
So, a sequence of flips terminates by the induction on $e(X,B)$. 
\end{step}
We complete the proof of Theorem \ref{main}. 
\end{proof}

\begin{rem}[Backtracking Method]\label{exist}
It is not difficult to see that the existence of 
$4$-dimensional semi-stable terminal flips implies that of 
all the $4$-dimensional semi-stable log flips. 
It is essentially proved in the proof of Theorem \ref{main}. 
We leave the details to the readers. 
See \cite[6.4, 6.5, 6.11 Theorem]{FA}. 
\end{rem}

\ifx\undefined\bysame
\newcommand{\bysame|{leavemode\hbox to3em{\hrulefill}\,}
\fi

\end{document}